\documentclass[12pt]{article}
\usepackage{graphicx}
\usepackage{amsmath}
\usepackage{amsfonts}
\usepackage{amssymb}

\newtheorem{theorem}{Theorem}

\newtheorem{lemma}[theorem]{Lemma}

\newtheorem{proposition}[theorem]{Proposition}

\newenvironment{proof}[1][Proof]{\textbf{#1.} }{\ \rule{0.5em}{0.5em}}

\def\text{\hbox} 
\date{}
\def\tr{{\rm tr\,}}
\def\a{\alpha}
\def\b{\beta}
\def\g{\gamma}

\def\d{\delta}

\def\l{\lambda}

\def\ro{\rho}

\def\ps{\varphi}

\def\G{\Gamma}

\def\A{{\cal A}}

\def\A'{{\cal A}'}

\def\wti{\widetilde}

\begin{document}

\title{Volumes for twist link cone-manifolds}

\author{D. Derevnin\footnote{Partially supported by RFBR (grant 03-01-00104),
INTAS (grant 03-51-3663) and Scientific Schools (grant
SS-300.2003.1).}
 \addtocounter{footnote}{-1} \and
A. Mednykh\footnotemark \and M. Mulazzani}
%\date{ }

 \maketitle

%\centerline{(preliminary version)}

\begin{abstract}
Recently, the explicit volume formulae for hyperbolic
cone-manifolds, whose underlying space is the 3-sphere and the
singular set is the knot $4_1$ and the links $5^2_1$ and $6^2_2$,
have been obtained by the second named author and his
collaborators. In this paper we explicitly find the hyperbolic
volume for cone-manifolds with the link $6^2_3$ as singular set.
Trigonometric identities (Tangent, Sine and Cosine Rules) between
complex lengths of singular components and cone angles are
obtained for an infinite family of two-bridge links containing
$5^2_1$ and $6^2_3$.
\\ {{\it Mathematics Subject Classification 2000:} Primary 57M50; Secondary 57M25, 57M27.}

\end{abstract}

\section{Introduction}

Starting from Alexander's works, polynomial invariants have become
a very convenient instrument for knot investigation. Several kinds
of knots polynomials have been discovered in the last twenty
years. Among these, we recall the
Jones-,~Kaufmann-,~HOMFLY-,~A-polynomials and others (\cite{Kauf},
\cite{CCGLS}, \cite{HLM2}). These polynomials relate knot theory
to algebra and algebraic geometry. Algebraic techniques are used
to find the most important geometrical characteristics of knots,
such as volume, length of shortest geodesics and others.

The explicit volume formulae for hyperbolic cone-manifolds, whose
underlying space is the 3-sphere and the singular set is the knot
$4_1$ and the links $5^2_1$ and $6^2_2$, have been obtained in
\cite{MR}, \cite{MV2} and \cite{M2}.

The aim of our paper is to explicitly find the hyperbolic volume
for cone-manifolds with the link $6^2_3$ as singular set. In order
to do this, we will introduce a family of hyperbolic
cone-manifolds $W_p(\a,\b)$, with the two-bridge links $W_p$, with
slope $(4p+4)/(2p+1)$ as singular set, and $\a,\b$ as cone angles.

Trigonometric identities (Tangent, Sine and Cosine Rules) between
complex lengths of singular components and cone angles for
$W_p(\a,\b)$ are obtained. Then the Schl\"afli formula applies in
order to find explicit hyperbolic volumes for cone-manifolds
$W_2(\a,\b)$.

In the present paper links and knots are considered  as singular
subsets of the three-sphere endowed by a Riemannian metric of
negative constant curvature.

\section{Trigonometric identities for knots and links}

\subsection{Cone-manifolds, complex distances and lengths}

We start with the definition of cone-manifold modelled in
hyperbolic, spherical or Euclidian structure.

\medskip

\noindent {\bf Definition 1}. A 3-dimensional {\it hyperbolic
cone-manifold} is a Riemannian 3-dimensional manifold of constant
negative sectional curvature with cone-type singularity along
simple closed geodesics.

\medskip

To each component of the singular set is associated a real number
$n \ge 1$  such that the cone-angle around the component is
$\alpha=2\pi/n .$ The concept of hyperbolic cone-manifold
generalizes that of hyperbolic manifold, which appears in the
partial case when all cone-angles are $2\pi.$ Hyperbolic
cone-manifolds are also a generalization of hyperbolic
3-orbifolds, which arises when all associated numbers $n$ are
integers. Euclidean and spherical cone-manifolds are defined
similarly.

In the present paper  hyperbolic, spherical or Euclidean
cone-manifolds $C$ are considered whose underlying space is the
three-dimensional sphere and the singular set $\Sigma = \Sigma^1
\cup \Sigma^2\cup \ldots \cup \Sigma^k $ is a link consisting of
the components $\Sigma^j = \Sigma^j(\alpha_j),\, j=1,2,\ldots,k$
with cone-angles $\alpha_1, \ldots, \alpha_k$  respectively.

We recall a few well-known facts from hyperbolic geometry.
\medskip

Let $\mathbf H^3=\{(z,\xi)\in \mathbf C\times \mathbf R : \xi>0
\}$ be the upper half space model of the 3~-di\-men\-sio\-nal
hy\-per\-bolic space en\-dowed by the Rieman\-nian metric
$$ds^2=\frac{dz d \overline z + d\xi^2}{\xi^2}.$$ We identify the
group of orien\-tation preserving isometries of $\mathbf H^3$ with
the group $PSL(2,\mathbf C)$, consisting of linear fractional
transformations $$        A': z\in \mathbf C \to \frac{az+b}{cz+d}
.$$ By a canonical procedure, $A'$ can be uniquely extended to an
isometry of $\mathbf H^3$. We prefer to deal with the matrix
$A=\left(
\begin{array}{cc}
a &  b \\
c & d
\end{array}
\right) \in SL(2,\mathbf C) $ \ rather than the element $A'\in
PSL(2,\mathbf C)$. The matrix $A$ is uniquely determined by the
element $A'$, up to a sign. In the following we will use the same
letter $A$ for both $A$ and $A'$, as long as this does not create
confusion.

Let $C$  be a hyperbolic cone-manifold with the singular set
$\Sigma.$ Then $C$ defines a nonsingular but incomplete hyperbolic
manifold ${\cal M} = C - \Sigma.$ Denote by $\Phi$ the fundamental
group of the manifold ${\cal M}$.

The hyperbolic structure of ${\cal M}$ defines, up to conjugation
in $PSL(2,\mathbf C)$, a holonomy homomorphism $$ \hat h  : \Phi
\to PSL(2,\mathbf C). $$ It is shown in \cite{Zhou} that the
holonomy homomorphism of an orientable cone-manifold can be lifted
to $SL(2,\mathbf C)$  if all cone-angles are at most $\pi$. Denote
by $h : \Phi \to SL(2,\mathbf C)$ this lifting homomorphism.
Choose an orientation on the link $\Sigma = \Sigma^1 \cup
\Sigma^2\cup \ldots \cup \Sigma^k $ and fix a meridian-longitude
pair $\{m_j,\,l_j\}$ for each component
$\Sigma^j=\Sigma^j(\alpha_j).$ Then the matrices $M_j=h(m_j)$ and
$L_j=h(l_j)$ satisfy the following properties: $$
\tr(M_j)=2\cos(\alpha_j/2),\,\,\,\, M_j L_j=L_j M_j\,,
\,\,\,j=1,2,\ldots,k. $$

Now we point out some definitions and results from the book
\cite{Fench}. A matrix $A \in SL(2,\mathbf C)$ satisfying  $\tr(A)
= 0$ is called a (normalized) {\it line matrix}. We have from
definition $A^2 = -I,$ where $I$ is the identity matrix. Hence any
line matrix determines a half-turn about a line in $\mathbf H^3,$
and this line determines the matrix up to sign. According to
\cite[p. 63]{Fench}, there exists a natural one-to-one
correspondence between line matrices and oriented lines in
$\mathbf H^3.$ Hereby, if a line matrix $A$ determines an oriented
line $\lambda_A=[e,e']$ with end points $e$ and $e'$, then the
line matrix $-A$ determines the line $[e',e]$. Moreover, if a
matrix $F \in SL(2,\mathbf C)$ is considered as a motion of
$\mathbf H^3$, then the matrix $FAF^{-1}$ determines the line
$[F(e),F(e')].$

\medskip

\noindent {\bf Definition 2}. Let $\lambda_A$ and $\lambda_B$ be
oriented lines determined by the line matrices $A$ and $B$. A
complex number $\mu$ is called a {\it complex distance} from
$\lambda_A$ to $\lambda_B$ if its real part $\Re\, \mu$ is the
distance from $\lambda_A$ to $\lambda_B$, and its imaginary part
$\Im\, \mu$ is the angle from $\lambda_A$ to $\lambda_B$ chosen in
$[0,2 \pi)$ .

\medskip

We have \cite[p. 68]{Fench}
\begin{equation} \label{ro}
  \cosh
\mu = -\frac 12\, \tr(AB).
\end{equation}
{}From now on, all lines in this paper will be supposed to be
oriented.

Any isometry $A$ of $\mathbf H^3$ different from parabolic and the
identity has two fixed points $u$ and $v$ in $\widehat{\mathbf
C}$. It acts as a translation of distance $r_A$ along the axis
$\l_A=[u,v]$ and rotation of $\ps_A$ about $\l_A$.

\medskip

\noindent {\bf Definition 3}. We call {\it displacement} of $A$
the complex number $\d(A)=r_A+i\ps_A$.

\medskip

The isometry $A$, without an orientation of its axis, determines
$\delta(A)$ up to sign. By \cite[p. 46]{Fench}, for the isometry
given by a matrix $A \in SL(2,\mathbf C)$ we have $$
 2 \cosh {\delta({A})}=\tr(A^2)=\tr^2(A)-2 .
$$

We remark that if $\d(A)\neq 0$ then $A$ has two different fixed points, so it admits
an axis determined by these points. The line matrix $\wti A$ of this axis is defined
by
\begin{equation} \label{three}
\wti A=\frac{A-A^{-1}}{2i\sinh \frac{\d(A)}{2}}  \end{equation}
Since $\d(A^{-1})=-\d(A),$ the matrices $A$ and $A^{-1}$ define
the same line matrix $\wti A=\wti {A^{-1}}~$  (see \cite{Fench}).

\medskip

\noindent {\bf Definition 4}. The {\it complex length} $\gamma_j$
of a singular component $\Sigma^j$ of the cone-manifold $C$ is the
displacement $\d(L_j)$ of the isometry $L_j$, where $L_j=h(l_j)$
is represented by the longitude $l_j$ of $\Sigma^j$.

\medskip

Immediately from the definition we get  \cite[p. 46]{Fench}
\begin{equation} \label{four}
2 \cosh \gamma_j= \tr (L_j^2).
\end{equation}

We note \cite[p. 38]{BZie} that the meridian-longitude pair $
\{m_j,\,l_j \}$ of the oriented link is uniquely determined up to
a common conjugating element of the group $\Phi$. Hence, the
complex length $\gamma_j=r_j+i \ \varphi_j$  is uniquely
determined $(\bmod \ 2\pi i)$, up to a sign, by the above
definition.

We need two conventions to correctly choose real and imaginary
parts of $\gamma_j.$ The first convention is the following. By the
assumptions on the singular set we have $r_j \ne 0.$ Hence, we can
choose $\gamma_j$ in such a way that $r_j>0.$ The second
convention concerns the imaginary part $\varphi_j.$ We want to
choose $\varphi_j $ such that the following identity holds
\begin{equation} \label{five}
\cosh \frac{\gamma_j}2 = -\frac 12 \,\tr (L_j)  \end{equation}

By virtue of identity ${\tr}^2 (L_j) - 2 = \tr (L_j^2)$ equality
(\ref{four}) is a consequence of (\ref{five}), but the converse,
in general, is true only up to a sign. Under the second convention
(\ref{four}) and (\ref{five}) are equivalent. The two above
conventions lead to convenient analytic formulas in order to
calculate $\g_j$ and $r_j$. More precisely, there are simple
relations between these numbers and the eigenvalues of the matrix
$L_j$. Recall that $\det (L_j) =1.$ Since $L_j$ is loxodromic, it
has two eigenvalues $f_j$ and $1/f_j$. We choose $f_j$ so that
$|f_j|>1.$ The case $|f_j|=1$ is impossible because in this case
the matrix $L_j$ is elliptic and therefore $r_j=0.$ Hence $$
f_j=-e^{\frac{\gamma_j}2}, \,\,\, |f_j|=e^{\frac{r_j}2}. $$

\bigskip
\begin{figure}[htbp]
\begin{center}
\includegraphics*[totalheight=7cm]{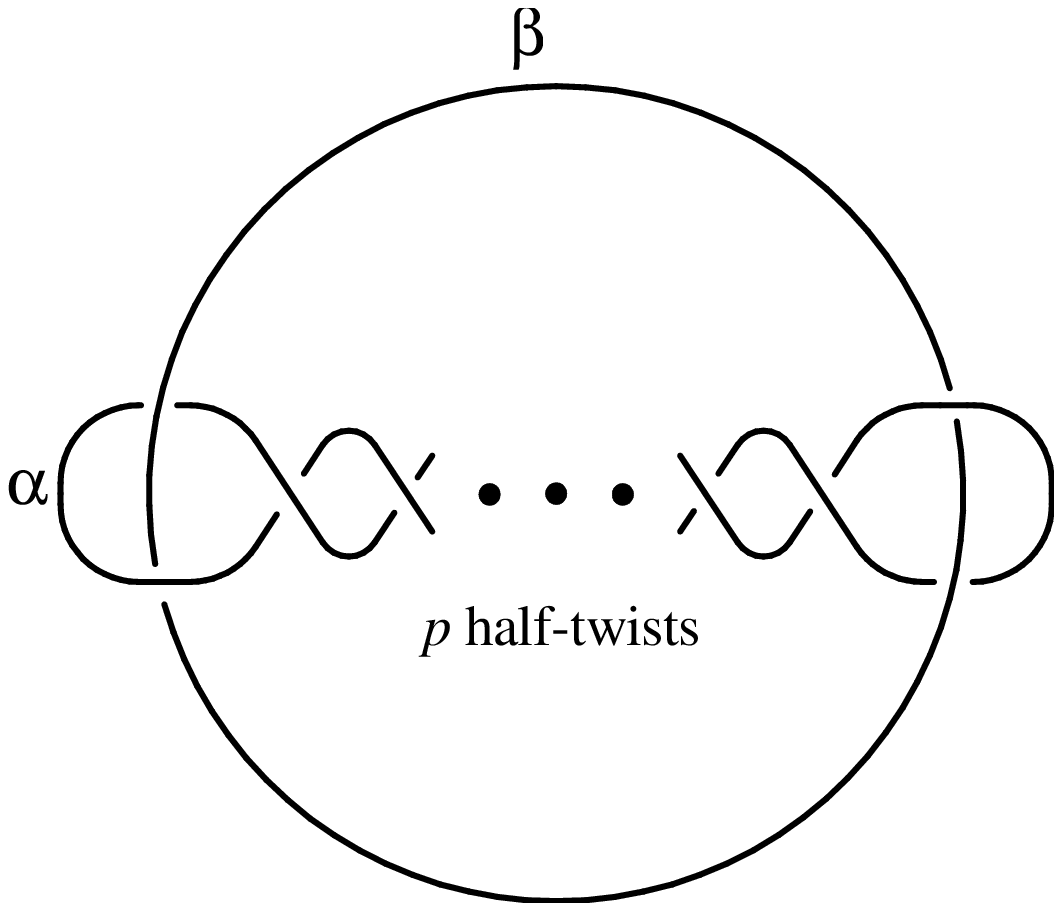}
\end{center}
\caption{The cone-manifold $W_p(\a,\b)$.} \label{Fig. 1}
\end{figure}

In this paper we consider a family of cone-manifolds whose
singular sets are links which are generalizations of the Whitehead
link. The link $W_p$, $p\ge 0$, is the two-component link depicted
in Figure~1, where $p$ is the number of half twists of one
component. For this reason we will call them {\it twist links}. It
is easy to see that $W_0$ is the torus link of type $(2,4)$ and
$W_1$ is the Whitehead link. All twist links are two-bridge links,
in particular $W_p$ is the two-bridge link with slope
$(4p+4)/(2p+1)$, for all $p\ge 0$. They are all hyperbolic, except
for $W_0$.

Denote by $W_p(\alpha,\beta)$ the cone-manifold whose underlying
space is the 3-sphere and whose singular set consists of the twist
link $W_p$ with cone angles $\alpha=2\pi/ m$ and $\beta=2\pi/n$
(see  Figure 1). It follows from the Thurston theorem that
$W_p(\alpha,\beta)$, with $p\neq 0$, admits a hyperbolic structure
for all sufficiently small $\alpha$ and $\beta$.

By the Kojima rigidity theorem \cite{Kj} the hyperbolic structure
is unique, up to isometry, if $0\le\a,\b\le\pi$.

In our paper we deal only with this range of angles.

Let us investigate the hyperbolic structure of the cone-manifold
$W_p(\a,\b)$. Its singular set $\Sigma = \Sigma^1 \cup \Sigma^2$
of consists of two components $\Sigma^1 = \Sigma^1(\a)$ and
$\Sigma^2 = \Sigma^2(\b)$ with cone-angles $\a$ and $\b$
respectively. $W_p(\a,\b)$ defines a nonsingular but incomplete
hyperbolic manifold ${\cal M} = W_p(\a,\b) - \Sigma.$ The
fundamental group of the manifold ${\cal M}$ has the following
presentation $$ \Phi_p=\langle s,t\mid\ sl_s=l_ss\rangle = \langle
s,t\mid\ tl_t=l_tt\rangle , $$ where $s$ and $t$ (resp. $l_s$ and
$l_t$) are meridians  (resp. longitudes) of the components
$\Sigma^1$ and $\Sigma^2$ respectively.

We use the following expression of $l_s$ in terms of $s$ and $t$:

\begin{equation} \label{odd}
l_s=[s,t]^{\frac{p+1}{2}}[s,t^{-1}]^{\frac{p+1}{2}}, \quad{\rm if }\,\,\,p\,\,\,{\rm
is\,\, odd}, \end{equation}
\begin{equation} \label{even}
l_s=s^{-1}[t,s]^{\frac p2}tst[s^{-1},t^{-1}]^{\frac p2}, \quad{\rm if
}\,\,\,p\,\,\,{\rm is\,\, even},\end{equation} where $[s,t]=sts^{-1}t^{-1}.$

The expressions for $l_t$ can be easily obtained by exchanging $s$ and $t$ in the
previous formulae.

Let $$\hat h=\hat h_{\a,\b} : \Phi_p \to PSL(2,\mathbf C) $$ and
$$h=h_{\a,\b}: \Phi_p \to SL(2,\mathbf C)$$ be holonomy
homomorphisms and $\Gamma_{\a,\b}= h_{\a,\b}(\Phi_p)$. The images
$\hat h_{\a,\b}(s)$ and $\hat h_{\a,\b}(t)$ of $s$ and $t$ are
rotations in $\mathbf H^3$ of angles $\a$ and $\b$ respectively.
The group $\Gamma_{\a,\b}$ is generated by the two matrices
$S=h_{\a,\b}(s)$ and $T=h_{\a,\b}(t)$ with the following
properties: $$ \tr(S)=2\cos \frac{\a}2,\quad \tr(T)=2\cos
\frac{\b}2,\quad \ SL_S=L_SS, $$ where $L_S=h_{\a,\b}(l_s).$

\subsection{Complex distance equation for two-bridge links}

The fundamental group of (the exterior of) a link $K$ is generated
by two meridians if and only if $K$ is a two-bridge link
\cite{BZ}. Moreover, a two-bridge link is hyperbolic if and only
if its slope is different from $p/1$ and $p/(p-1)$ (see \cite{T}).

\begin{proposition} \label{equation} Let
$ \Phi=\langle s,t\rangle$ be the fundamental group of a
hyperbolic two-bridge link $K$ generated by the two meridians $s$
and $t$. Let $\G_{\a,\b}=h_{\a,\b}(\Phi)$ be the image of $\Phi$
under the holonomy homomorphism of the hyperbolic cone manifold
$K(\a,\b)$. Then, up to conjugation in $SL(2,\mathbf C),$ the
generators $S=h_{\a,\b}(s)$ and $T=h_{\a,\b}(t)$ of $\G_{\a,\b}$
can be chosen in such a way that
\begin{equation} \label{six}
S=\left(
\begin{array}{cc}
\cos\frac{\a}2 & i\,e^{\frac{\rho}{2}} \sin \frac{\a}2  \\
i\,e^{-\frac{\rho}{2}} \sin \frac{\a}2  & \cos \frac{\a}2 \\
\end{array}
\right)\ ,\ \ T=\left(
\begin{array}{cc}
\cos \frac{\b}2 & i\,e^{-\frac{\rho}{2}} \sin\frac{\b}2  \\
i\,e^{\frac{\rho}{2}} \sin \frac{\b}2  & \cos \frac{\b}2 \\
\end{array}
\right),
\end{equation}
where $\ro$ is the complex distance between the axis of $S$ and
$T$.
\end{proposition}
\begin{proof}  After a suitable conjugation in the group $SL(2,\mathbf
C),$ one can assume that the oriented axes of the elliptic
elements $S$  and $T$ are $\l_S=[-e^{\frac{\ro}2},
e^{\frac{\ro}2}]$ and $\l_T=[-e^{-\frac{\ro}2},e^{-\frac{\ro}2}]$.
Since $\tr(S)=2\,\cos \frac{\a}2$ \ and \ $\tr(T)=2\,\cos
\frac{\b}2$, \, the matrices $S$ and $T$ are given by (\ref{six}).
Check that $\rho$ coincides with the complex distance $\rho(S,T)$
between $\l_S$ and $\l_T$. The line matrices $\wti S$ and $\wti
T$, corresponding to these axes, can be obtained by (\ref{three}).
Since $\d(S)=i\,\a$ and $\d(T)=i\,\b$, we have $\wti S=\left(
\begin{array}{cc}
0 & -i e^{\frac{\rho}{2}} \\
-i e^{-\frac{\rho}{2}} & 0 \\
\end{array}
\right) \   $ and $\wti T=\left(
\begin{array}{cc}
0 & -i e^{-\frac{\rho}{2}} \\
-i e^{\frac{\rho}{2}} & 0 \\
\end{array}
\right)$ respectively.  By \cite[p. 68]{Fench} we get $\cosh
\ro(S,T)=-\frac 12 \ \tr(\wti S\wti T)=\cosh \ro$.
\end{proof}

\medskip

The following two propositions can be obtained by direct
calculation from the above statement.
\begin{proposition}\label{equation1}Let $$ \Phi_2=\langle s,t: sl=ls,
\,l=s^{-1}tst^{-1}s^{-1}tsts^{-1}t^{-1}st\rangle$$ be the
fundamental group of the two-bridge link $W_2$ with slope $12/5$
and $\G_{\a,\b}=h_{\a,\b}(\Phi_2)=\langle S,T\rangle$ be the image
of $\Phi_2$ under the holonomy homomorphism of the hyperbolic cone
manifold $W_2(\a,\b)$. Denote by $\ro= \ro(S,T)$ the complex
distance between the axes of $S=h_{\a,\b}(s)$ and
$T=h_{\a,\b}(t).$ Then $u=\cosh \ro$ is a non-real root of the
complex distance equation
\begin{equation} \label{phi2}
4z^3 - 4ab z^2 + (3a^2b^2 + 3a^2 + 3b^2 -1) z - ab(a^2b^2 + a^2 + b^2 -3) =0,
\end{equation}
 where $a=\cot \frac{\a}2$ and $b=\cot \frac{\b}2$.
\end{proposition}
\begin{proof}
Denote by $ L=S^{-1}TST^{-1}S^{-1}TSTS^{-1}T^{-1}ST$ the image of
the longitude $l$ under the holonomy homomorphism $h = h_{\a,\b}:
\Phi_2 \to SL(2,\mathbf C).$  Then we have $SL=LS$.

Let $N$ be a line matrix corresponding to the common normal to the
axes of $S$ and $T$. If $S$ and $T$ are represented in the form
(\ref{six}) then one can take $N=\left(
\begin{array}{cc}
i & 0 \\
0 & -i \\
\end{array}
\right) \ . $ It is not difficult to verify that $NSN^{-1}=S^{-1}$ and
$NTN^{-1}=T^{-1}.$

To complete the proof, we need the following lemma, which gives
simple criteria for matrices $S$ and $L$ to be permutable.
\begin{lemma} The following conditions are equivalent: (i) $SL=LS;$\,
(ii) $NLN^{-1}=L^{-1};$ \, (iii) $\tr(NL)=0.$
\end{lemma}
\begin{proof}  First we show that (i) and (ii) are
equivalent. Indeed, since $ L=S^{-1}TST^{-1}S^{-1}TSTS^{-1}T^{-1}ST$ we have
$NLN^{-1}=ST^{-1}S^{-1}TST^{-1}S^{-1}T^{-1}STS^{-1}T^{-1}=SL^{-1}S^{-1}.$ Hence (ii)
holds if and only if $S$ and $L^{-1}$ are permutable. The last property is equivalent
to (i). Because of $N^2=-I$ the condition (ii) can be rewritten in the form $NLNL=-I;$
that is equivalent to (iii).
\end{proof}

By this lemma and direct calculation we have $$ \tr(NL)= \frac{-4i\sinh
\ro}{(1+a^2)^3(1+b^2)^3}\,\cdot (4u^2+a^2b^2+a^2+b^2-3)\,\cdot $$ $$ \cdot\,(4u^3 -
4ab u^2 + (3a^2b^2 + 3a^2 + 3b^2 -1) u - ab(a^2b^2 + a^2 + b^2 -3)) =0, $$ where
$u=\cosh\rho$.

Now we have to show that $u$ is a non-real root of (\ref{phi2}).
Since $\G_{\a,\b}$ is the holonomy group of a hyperbolic
cone-manifold, it is non-elementary\footnote{A subgroup $G$ of
$SL(2,\mathbf C)$ is called {\it elementary} if it has a finite
orbit in ${\mathbf H}^3 \cup \widehat{\mathbf C}.$} and is not
conjugated to a subgroup of $SL(2,\mathbf R)$ \cite{HLM2}.

If $\sinh \ro = 0$ then the axes $S$ and $T$ coincide, and the group $\G_{\a,\b}$ is
elementary.

If $u$ is a root of equation $$4u^2+a^2b^2+a^2+b^2-3=0$$ then by
equality
$$\tr{L}-2=-\frac{4(a^2+u^2)(4u^2+a^2b^2+a^2+b^2-3)^2}{(a^2+1)^3(b^2+1)^3}
$$  we have $\tr{L}=2 .$ {}From (\ref{five}) we obtain $$\cosh
\frac{\gamma_S}2 = -\frac 12 \,\tr (L)=-1. $$ Hence
$\gamma_S=r_S+i\varphi_S=2\pi i$ and the real length $r_S$ of the
link component
 $\Sigma_1$ is equal to zero, which is a contradiction.

Suppose that $u=\cosh \ro$ is a real root. Let $$
R(z_1,z_2,z_3,z_4)=\frac{(z_3-z_1)(z_4-z_2)}{(z_3-z_2)(z_4-z_1)}$$ be the cross ratio
of the four points $z_1,z_2,z_3,z_4\in\widehat{\mathbf C}$. Then $R(-e^{\frac{\ro}2},
e^{\frac{\ro}2},-e^{-\frac{\ro}2},e^{-\frac{\ro}2}) = (\cosh \ro -1)/(\cosh \ro +1)
\in \mathbf R \cup \{\infty \}.$ We have that the axes $[-
e^{\frac{\ro}2},e^{\frac{\ro}2}]$ and $[-e^{-\frac{\ro}2},e^{-\frac{\ro}2}]$ of $S$
and $T$ lie in a common plane. If the axes intersect then  the group
$\G_{\a,\b}=\langle S,T \rangle$ has a fixed point and is elementary. If they do not
intersect, $\G_{\a,\b}$ is conjugated to a subgroup of $SL(2,\mathbf R)$.

Therefore, we have shown that $u$ is a non-real root of (\ref{phi2}) and the proof of
Proposition \ref{equation1} is completed.
\end{proof}

\medskip

The next proposition can be proved by similar arguments.
\begin{proposition}\label{equation2}Let $$\Phi_3=\langle s,t:~sl~=~ls,
\,~l~=~sts^{-1}t^{-1}sts^{-1}t^{-1}st^{-1}s^{-1}tst^{-1}s^{-1}t\rangle$$
be the fundamental group of the two-bridge link $W_3$ with the
slope $16/7$ and $\G_{\a,\b}=h_{\a,\b}(\Phi_3)=\langle S,T\rangle$
the image of $\Phi_3$ under the holonomy homomorphism of a
hyperbolic cone manifold $W_3(\a,\b)$ generated by
$S=h_{\a,\b}(s)$ and $T=h_{\a,\b}(t).$ Denote by $\ro= \ro(S,T)$
the complex distance between the axes of $S$ and $T.$ Then
$u=\cosh \ro$ is a non-real root of the complex distance equation
$$ 0=8u^5 + 8ab u^4 +8 (a^2b^2 + a^2 + b^2 -1) u^3 + 4ab(a^2b^2 +
a^2 + b^2 -3) u^2 +$$
$$(a^4b^4+2a^4b^2+2a^2b^4-4a^2b^2+a^4+b^4-6a^2-6b^2+1) u
-4ab(a^2b^2+a^2+b^2-1) , $$ where $a=\cot \frac{\a}2$ and $b=\cot
\frac{\b}2$.
\end{proposition}

\subsection{Tangent, Sine and Cosine rules}

If we set $z=\tr(S^{-1}T)$ then, from presentation in Proposition
\ref{equation}, we have
$$z=2(\cos\frac{\a}2\cos\frac{\b}2+u\sin\frac{\a}2\sin\frac{\b}2
),$$ where $u=\cosh\rho$.

The algebraic equation for $z$ and its behaviour was considered in
a number of papers (see \cite{CCGLS}, \cite{GM}, \cite{HLM2} and
others) devoted to $PSL(2,\mathbf C)$ representation of
two-generator groups.

In general, the equation for $u$ (as well as for $z$) is very
complicated, even for twist links. In spite of this, since
$u=\cosh\rho$ has a very clear geometric sense, we are able to
produce some general results for twist links without calculating
$u$.

\begin{proposition} \label{pizza} Let $W_p(\a,\b)$ be a hyperbolic twist link cone-manifold.
Denote by $S=h_{\a,\b}(s)$ and $T=h_{\a,\b}(t)$ the images of the
generators of the group $\Phi_p=\langle s,t\mid sl_s=l_ss\rangle$
under the holonomy homomorphism $h_{\a,\b}:\Phi_p\to SL(2,\mathbf
C)$. Set $u=\cosh{\ro}$, where $\ro$ is the complex distance
between the axes of $S$ and $T$, such that $\Im \, u>0$. Moreover,
denote by $\gamma_{\a}$ and $\gamma_{\b}$ the complex lengths of
the singular components of $W_p(\a,\b)$ with cone-angles $\a$ and
$\b$ respectively. Then $$ u=i\,\cot \frac{\a}2 \coth
\frac{\gamma_{\b}}4=i\,\cot \frac{\b}2 \coth
\frac{\gamma_{\a}}4.$$
\end{proposition}
\begin{proof} To prove the statement we need to calculate the complex distance
between axes of elliptic elements $S$ and $T$ in two ways. By
definition, $L_S=h_{\a,\b}(l_s)$ and $L_T=h_{\a,\b}(l_t)$, where
$l_s$ and $l_t$ are the longitudes of the singular components of
$W_p(\a,\b)$ with cone angles $\a$ and $\b$ respectively.

First of all we fix an orientation on the axes of $S$ and $T$ by
the following line matrices $$\wti S=\frac{S-S^{-1}}{2\,i\,\sinh
\frac{i\,\a}{2}},\qquad \wti T =\frac{T-T^{-1}}{2\,i\,\sinh
\frac{i\,\b}{2}}.$$ Then the complex distance $\rho(S,T)$ between
the oriented axes of $S$ and $T$ is defined by (\ref{ro}):
$$\cosh\rho(S,T)=-\frac12\,\tr(\wti S\wti T).$$ Using
(\ref{three}) we define the line matrices for $L_S$ and $L_T$ as
$$\wti {L_S}=\frac{L_S-L_S^{-1}}{2i\sinh \frac{\g_{\a}}{2}},
\qquad \wti {L_T}=\frac{L_T-L_T^{-1}}{2i\sinh
\frac{\g_{\b}}{2}}.$$ To continue the proof, we need two lemmas:
\begin{lemma} \label{lemma1}
For every $S,T$ we have $\wti S=-\wti {L_S}$ and $\wti T=-\wti
{L_T}.$
\end{lemma}
\begin{proof}
Up to conjugation in $SL(2,\mathbf C)$, we can assume that $S$ is
given by $$ S=\left(
\begin{array}{cc}
e^{\frac{i\a}2} & 0 \\
0  & e^{-\frac{i\a}2} \\
\end{array}
\right) .$$ Note that $ L_S$ is a loxodromic element, with
displacement $\g_{\a}$, permutable with $S$. Since $ \wti
{L^{-1}_S}=\wti {L_S}~$, we can assume that $$ L_S=\left(
\begin{array}{cc}
\pm e^{\frac{\g_{\a}}2} & 0 \\
0  & \pm e^{-\frac{\g_{\a}}2} \\
\end{array}
\right) $$ By convention (see formula (\ref{five})) we have $$\tr (L_S)=-2\cosh
\frac{\g_{\a}}2.$$ Hence $$ L_S=\left(
\begin{array}{cc}
-e^{\frac{\g_{\a}}2} & 0 \\
0  & -e^{-\frac{\g_{\a}}2} \\
\end{array}
\right) $$ and we obtain $$\wti L_S=\frac{L_S-L_S^{-1}}{2i\sinh\frac{\g_{\a}}2}=\left(
\begin{array}{cc}
i & 0 \\
0  & -i \\
\end{array}
\right) $$ and $$\wti S=\frac{S-S^{-1}}{2\,i\sinh\frac{i\,\a}2}=\left(
\begin{array}{cc}
-i & 0 \\
0  & i \\
\end{array}
\right).
$$
\end{proof}

\begin{lemma} \label{lemma2}
For every $S,T$ we have $\tr(S)=\tr(S^{-1}L_T)$ and
$\tr(T)=\tr(T^{-1}L_S)$.
\end{lemma}
\begin{proof}
To prove $\tr(T)=\tr(T^{-1}L_S)$ it is enough to show that
$T^{-1}L_S$ is conjugated to either $T$ or $T^{-1}$ in the group
$\Gamma_{\a,\b}$. If $p$ is odd, we have from (\ref{odd}):
$$T^{-1}L_S=T^{-1}[S,T]^{\frac{p+1}{2}}\,[S,T^{-1}]^{\frac{p+1}{2}}=
[T^{-1},S]^{\frac{p+1}{2}}\,T^{-1}[T^{-1},S]^{-\frac{p+1}{2}}.$$
If $p$ is even, we have from (\ref{even}):
$$T^{-1}L_S=T^{-1}S^{-1}[T,S]^{\frac{p}{2}}\,TST\,[S^{-1},T^{-1}]^{\frac{p}{2}}=
T^{-1}S^{-1}[T,S]^{\frac{p}{2}}\,T\,[T,S]^{-\frac{p}{2}}\,ST.$$ The equality
$\tr(S)=\tr(S^{-1}L_T)$ can be obtained in a similar way.
\end{proof}

\medskip

To complete the proof of Proposition \ref{pizza}, we note that
$\tr(XY)=\tr(X)\tr(Y)-\tr(X^{-1}Y)$, $\tr(X^{-1})=\tr(X)$ and
$\tr(XY)=\tr(X^{-1}Y^{-1})$ holds for all $X,Y\in SL(2,\mathbf
C)$. By Lemma \ref{lemma1}, Lemma \ref{lemma2} and formulae
$\tr(S)=2\cos\frac{\a}2$, $\tr (L_S)=-2\cosh \frac{\g_{\a}}2$, we
have $$\cosh\rho(S,T)= -\frac12\,\tr(\wti S \wti
T)=\frac12\,\tr(\wti S\wti {L_T})=$$
$$=\frac12\,\tr\left(\frac{S-S^{-1}}{2\sin \frac{\a}{2}}\,
\frac{L_T-L_T^{-1}}{2i\sinh \frac{\g_{\b}}{2}}\right)
=\frac{\tr(SL_T-S^{-1}L_T-SL_T^{-1}+S^{-1}L_T^{-1})}{8i\sin\frac{\a}2\sinh\frac{\g_{\b}}2}=$$
$$=\frac{2(\tr(SL_T)-\tr(S^{-1}L_T))}{8i\sin\frac{\a}2\sinh\frac{\g_{\b}}2}
=\frac{\tr(S)\tr(L_T)-2\tr(S^{-1}L_T)}{4i\sin\frac{\a}2\sinh\frac{\g_{\b}}2}=$$
$$=\frac{\tr(S)\tr(L_T)-2\tr(S)}{4i\sin\frac{\a}2\sinh\frac{\g_{\b}}2}=
\frac{\tr(S)(2-\tr(L_T))}{-4i\sin\frac{\a}2\sinh\frac{\g_{\b}}2}$$
$$=\frac{2\cos\frac{\a}2(2+2\cosh
\frac{\g_{\b}}2)}{-4i\sin\frac{\a}2\sinh\frac{\g_{\b}}2}
=i\cot\frac{\a}2 \coth \frac{\g_{\b}}4.$$ Since $\cosh
\rho(S,T)=\cosh \rho(T,S)=u$ the statement follows.
\end{proof}

\medskip

As an immediate consequence of the previous proposition, we have
the following result.

\begin{theorem} (The Tangent Rule) Suppose that $W_p(\alpha,\beta)$ is a
hyperbolic cone-manifold. Denote by $\gamma_\alpha$ and
$\gamma_\beta$ complex lengths of the singular geodesics of
$W_p(\alpha,\beta)$ with cone angles $\alpha$ and $\beta$
respectively. Then $$ \frac{\tanh{\frac{\gamma_\alpha}4}}
{\tanh{\frac{\gamma_\beta}4}} = \frac{\tan{\frac\alpha
2}}{\tan{\frac\beta 2 }}. $$
\end{theorem}

The following two theorems are consequences of the Tangent Rule.

\begin{theorem} (The Sine Rule)
Let $\gamma_\alpha =r_\alpha +i\,\varphi_\alpha $ and $\gamma_\b =r_\b +i\,\varphi_\b$
be the complex lengths of the singular geodesics of a hyperbolic cone-manifold
$W_p(\alpha,\beta)$ with cone angle $\alpha$ and $\beta$ respectively. Then $$
\frac{\sin{\frac{\varphi_\alpha}2}}{\sinh{\frac{r_\alpha}2}} =
\frac{\sin{\frac{\varphi_\beta}2}}{\sinh{\frac{r_\beta}2}}. $$
\end{theorem}

\begin{proof}By the Tangent Rule we have
$$ \frac{\tanh{\frac{\gamma_\alpha}4}}{a} = \frac{\tanh{\frac{\gamma_\beta}4}}{b}, $$
where $a=\tan{\frac{\alpha}2}$ and $B=\tan{\frac{\beta}2}$ are real numbers. Hence $$
\frac{\Re(\tanh{\frac{\gamma_\alpha}4})}{a} =
\frac{\Re(\tanh{\frac{\gamma_\beta}4})}{b}, $$ and $$
\frac{\Im(\tanh{\frac{\gamma_\alpha}4})}{a} =
\frac{\Im(\tanh{\frac{\gamma_\beta}4})}{b}. $$ Dividing one equation by the other we
obtain $$ \frac{\Re(\tanh{\frac{\gamma_\alpha}4})}{\Im(\tanh{\frac{\gamma_\alpha}4})}
=
\frac{\Re(\tanh{\frac{\gamma_\beta}4})}{\Im(\tanh{\frac{\gamma_\beta}4})}.
$$

By direct calculations we have
$$ \Re(\tanh{\frac{\gamma_\alpha}4})=\frac12
(\tanh{\frac{\gamma_\alpha}4}+
\tanh{\frac{\bar{\gamma}_\alpha}4})
%=\frac{\sinh(\frac{\gamma_\alpha}4+\frac{\bar{\gamma}_\alpha}4)}
%{2\,\cosh{\frac{\gamma_\alpha}4}\cosh{\frac{\bar{\gamma}_\alpha}4}}
=\frac{\sinh{\frac{r_\alpha}2}}
{\cosh{\frac{r_\alpha}2}+\cos{\frac{\varphi_\alpha}2}} $$ and $$
\Im(\tanh{\frac{\gamma_\alpha}4})=\frac{1}{2i}
(\tanh{\frac{\gamma_\alpha}4}- \tanh{\frac{\bar{\gamma}_\alpha}4})
%=\frac{\sinh(\frac{\gamma_\alpha}4-\frac{\bar{\gamma}_\alpha}4)}
%
%{2\,\cosh{\frac{\gamma_\alpha}4}\cosh{\frac{\bar{\gamma}_\alpha}4}}
=\frac{\sin{\frac{\varphi_\alpha}2}}
{\cosh{\frac{r_\alpha}2}+\cos{\frac{\varphi_\alpha}2}}. $$ Since
$r_\alpha > 0,$ we have $\displaystyle{\cosh{\frac{r_\alpha}2} >
1.}$ Therefore
$\displaystyle{\cosh{\frac{r_\alpha}2}+\cos{\frac{\varphi_\alpha}2}
> 0}$ and the result follows.
\end{proof}

\medskip

\begin{theorem} (The Cosine Rule)
Let $\gamma_\alpha =r_\alpha +i\,\varphi_\alpha $ and $\gamma_\b =r_\b +i\,\varphi_\b$
be the complex lengths of the singular geodesics of a hyperbolic cone-manifold
$W_p(\alpha,\beta)$ with cone angle $\alpha$ and $\beta$ respectively. Then $$ \frac{
\cos{\frac{\varphi_\alpha}2}\cosh{\frac{r_\beta}2}
-\cos{\frac{\varphi_\beta}2}\cosh{\frac{r_\alpha}2}}
{\cosh{\frac{r_\alpha}2}\cosh{\frac{r_\beta}2}
-\cos{\frac{\varphi_\alpha}2}\cos{\frac{\varphi_\beta}2}} =\frac{\cos\a-\cos\b
}{1-\cos\a\cos\b }. $$
\end{theorem}
\begin{proof}
By the Tangent Rule $$
\frac{\tanh{\frac{\gamma_\alpha}4}\tanh{\frac{\bar{\gamma}_\alpha}4}}{a^2} =
\frac{\tanh{\frac{\gamma_\beta}4}\tanh{\frac{\bar{\gamma}_\beta }4}}{b^2}, $$ where
$a=\tan{\frac{\alpha}2}$ and $b=\tan{\frac{\beta}2}.$ Hence $$
\frac{1+\cos\a}{1-\cos\a}\,\frac{\cosh{\frac{r_\alpha}2}-\cos{\frac{\varphi_\alpha}2}}
{\cosh{\frac{r_\alpha}2}+\cos{\frac{\varphi_\alpha}2}}=\frac{1+\cos\b}{1-\cos\b}\,
\frac{\cosh{\frac{r_\beta}2}-\cos{\frac{\varphi_\beta}2}}
{\cosh{\frac{r_\beta}2}+\cos{\frac{\varphi_\beta}2}}. $$ Set $$p=\cos\a,\quad
q=\cos\b,\quad p'=\frac{\cos{\frac{\varphi_\alpha}2}}{\cosh{\frac{r_\alpha}2}},\quad
q'=\frac{\cos{\frac{\varphi_\beta}2}}{\cosh{\frac{r_\beta}2}}$$ and rewrite the above
equation in the form $$
\frac{1+p}{1-p}\,\frac{1-p'}{1+p'}=\frac{1+q}{1-q}\,\frac{1-q'}{1+q'}, $$ or,
equivalently, as $$
\log\frac{1+p}{1-p}+\log\frac{1-p'}{1+p'}=\log\frac{1+q}{1-q}+\log\frac{1-q'}{1+q'}.
$$ Since $\displaystyle{{\rm arctanh}\,p=\frac12\,\log\frac{1+p}{1-p}}$ we have $${\rm
arctanh}\,p-{\rm arctanh}\,p'={\rm arctanh}\,q-{\rm
arctanh}\,q'.$$ and $${\rm arctanh}\,p-{\rm arctanh}\,q={\rm
arctanh}\,p'-{\rm arctanh}\,q'.$$ Hence $$
\frac{p-q}{1-pq}=\frac{p'-q'}{1-p'q'} $$ and, after substituting
the expressions for $p,q,p',q'$ in the last formula, we obtain the
statement.
\end{proof}

\medskip

We remark that, in the case of Whitehead link cone-manifolds,
Tangent and Sine rules are obtained in \cite{M1}.

\section{Explicit volume calculation for twist link cone-manifolds}

\subsection{The Schl\"afli formula}

In this section we will obtain explicit formulae for the volume of
some special cone-manifolds in the hyperbolic and spherical
geometries. In the case of complete hyperbolic structure on the
simplest knot and link complements such formulas, in terms of
Lobachevsky function, are well-known and widely represented in
\cite{T}. In general, a hyperbolic cone-manifold can be obtained
by completion of a non-complete hyperbolic structure on a suitable
knot or link complement. If the cone-manifold is compact, explicit
formulas are only known in a few cases \cite{HLM3}, \cite{HLM4},
\cite{Hds}, \cite{M2}, \cite{M3}, \cite{MR}, \cite{MV1},
\cite{MV2}. In all these cases the starting point for the volume
calculation is the Schl\"afli formula (see, for example
\cite{Hds}).

\medskip

\begin{theorem} (The Schl\"afli volume formula) Suppose
that $C_t$ is a smooth 1-parameter family of (curvature $K$)
cone-manifold structures on an $n$-manifold, with singular locus
$\Sigma$ of a fixed topological type. Then the derivative of
volume of $C_t$ satisfies $$ (n-1)K dV(C_t) = \sum_{\sigma}
V_{n-2}(\sigma)\, d\theta(\sigma) $$ where the sum is over all the
components $\sigma$ of the singular locus $\Sigma,$ and
$\theta(\sigma)$ is the cone angle along $\sigma.$
\end{theorem}

In the present paper we will deal mostly with three-dimensional
cone-manifold structures of negative constant curvature $K=-1$.
The Schl\"afli formula in this case reduces to $$dV= -\frac12
\sum_i r_i d\theta_i,$$ where the sum is taken over all the
components of the singular set $\Sigma$ with lengths $r_i$ and
cone angles $\theta_i.$

Our aim is to obtain the volume formulas for twist link hyperbolic
cone-manifolds $W_2(\alpha, \beta)$. We note that the volume
formula for $W_1(\alpha, \beta)$ were obtained earlier in
\cite{M3} and \cite{MV2}.

\begin{proposition} \label{ultima}  Let $W_2(\alpha,\beta)$ be a hyperbolic
cone-manifold and $r_\alpha$, $r_\beta$ the lengths of its
singular components, with cone angles $\alpha$ and $\beta$
respectively. If $a=\cot\frac{\alpha}2$ and $b=\cot\frac{\beta}2$,
then
\begin{equation} \label{eight}
r_{\alpha}= 2i\,\arctan\frac{a}{\zeta}- 2i\,\arctan\frac{a}{\overline{\zeta}} ,
\end{equation}
\begin{equation} \label{nine}
 r_{\beta}= 2i\,\arctan\frac{b}{\zeta}-
2i\,\arctan\frac{b}{\overline{\zeta}} ,
\end{equation}
 where $\zeta$ is a root of the
equation
\begin{equation} \label{ten}
4(z^2 + a^2)(z^2 + b^2)-(1+a^2)(1+b^2)(z-z^2)^2=0,
\end{equation}
with $\Im(\zeta) > 0$.
\end{proposition}
\begin{proof}
By Proposition \ref{pizza} we have
\begin{equation} \label{eleven}
i\,b\, {\rm coth}\frac{\gamma_{\alpha}}{4}= i\,a\, {\rm
coth}\frac{\gamma_{\beta}}{4}=u, \end{equation} where
$u=\cosh{\ro}$, and $\ro$ is a complex distance between the axes
of $S$ and $T$, chosen so that $\Im \, u>0$. By Proposition
\ref{equation1}, $u$ is a root of the cubic equation $$ 4z^3 - 4ab
z^2 + (3a^2b^2 + 3a^2 + 3b^2 -1) z - ab(a^2b^2 + a^2 + b^2 -3) =0.
$$ {}From (\ref{eleven}), for a suitable choice of analytical
branches, $$ r_{\alpha}=\frac{\gamma_{\alpha}}2+
\frac{\overline{\gamma}_{\alpha}}2=
   2i\,\arctan\frac{\overline{u}}b- 2i\,\arctan\frac ub  =
2 i\,\arctan\frac{a}{\zeta}- 2i\,\arctan\frac{a}{\overline{\zeta}}
, $$ where $\zeta = ab/\overline{u}$, $\Im(\zeta) > 0$, satisfies
the equation
$$Q(z)=(a^2b^2+a^2+b^2-3)z^3 - (3a^2b^2 + 3a^2 + 3b^2 -1) z^2 + 4a^2b^2 z - 4a^2b^2
=0.$$ To finish the proof we note that $$(z+1)Q(z)=-4(z^2 + a^2)(z^2 +
b^2)+(1+a^2)(1+b^2)(z-z^2)^2.$$
\end{proof}

\medskip

In the next section we will apply this result to calculate the volume of
$W_2(\alpha,\beta)$ via the Schl\"afli formula.

We remark that formulae (\ref{eight}) and (\ref{nine}), as a
consequence of the Tangent Rule, also hold for all twist links
$W_p$, with $\zeta=ab/\bar u$, where $u=\cosh\rho$.

For example, an analog for the algebraic equation (\ref{ten}), in
the case of twist link $W_3$, can easily be obtained from
Proposition \ref{equation2}. But in this case the equation became
too complicated and we are not able to explicitly find the
integrand in the Schl\"afli formula.

\subsection{Volume of twist link cone-manifolds}

The case of the Whitehead link cone manifolds $W_1(\alpha, \beta)$ has already been
solved (see \cite{M3} and \cite{MV2}).

\begin{theorem} {\em \cite{M3,MV2}} Let $W_1( \alpha, \beta)$ be a hyperbolic
Whitehead link cone-manifold with cone angles $\alpha$  and
$\beta$. Then the volume of $W_1(\alpha, \beta)$ is given by the
formula $$ {\rm Vol}\,W_1( \alpha, \beta)=
i\int_{\overline{\zeta}}^{\zeta}\log \left[ \frac{2(z^2 + a^2)(z^2
+ b^2)}{(1+a^2)(1+b^2)(z^2-z^3)} \right]\frac{dz}{z^2-1}. $$ where
$a=\cot\frac{\alpha}2,\,b=\cot\frac{\beta}2$ and $\zeta$ is a
non-real root,  with $\Im(\zeta)>0$,  of the equation $$2(z^2 +
a^2)(z^2 + b^2)-(1+a^2)(1+b^2)(z^2-z^3)=0.$$
\end{theorem}

The main result of this section is the following.

\begin{theorem} Let $W_2( \alpha, \beta)$ be a hyperbolic
twist link cone-manifold with cone angles $\alpha$  and $ \beta$.
Then the volume of $W_2(\alpha, \beta)$ is given by the formula

\begin{equation} \label{more}
 {\rm Vol}\, W_2( \alpha, \beta)=i\,
\int_{\overline{\zeta}}^{\zeta}\log\left[ \frac{4(z^2 + a^2)(z^2 +
b^2)}{(1+a^2)(1+b^2)(z-z^2)^2}\right]\frac{dz}{z^2-1}. \end{equation}
where
$a=\cot\frac{\alpha}2$, $b=\cot\frac{\beta}2$ and  $\zeta$ is a non-real root, with
$\Im(\zeta)>0$, of the equation
\begin{equation} \label{prova}
4(z^2 + a^2)(z^2 + b^2)-(1+a^2)(1+b^2)(z-z^2)^2=0. \end{equation}
\end{theorem}

\begin{proof} Denote by $V = {\rm Vol}\, W_2( \alpha, \beta)$ the hyperbolic volume of
$W_2( \alpha, \beta).$ Then by virtue of the Schl\"afli formula we have
\begin{equation} \label{twelve}
\frac{\partial
V}{\partial \alpha}=-\frac{r_{\alpha}}2,\,\,\,
 \frac{\partial V}{\partial \beta}=-\frac{r_{\beta}}2,
\end{equation}
where $r_{\alpha}$ and $r_{\alpha}$ are the lengths of the
singular geodesics having cone angles $\alpha$ and $\beta$
respectively.

We note that  for $\alpha=\beta$  and $\Im(\zeta) \to 0$ the
geometrical limit of  the cone-manifold $W_2(\alpha,\,\alpha)$ is
an Euclidean cone manifold  $W_2(\alpha_0,\,\alpha_0),$  where
$\alpha_0=2.7243... <  \pi.$ (See Example 1 in Section 3.3 below).
Hence, by Theorem 7.1.2 of \cite{Kj}, we have
\begin{equation} \label{thirteen}
 V \to 0  \,\,\,{\rm as}\,\,\,  \alpha=\beta \,\, \,{\rm and}\,\,\,  \Im(\zeta) \to 0.
\end{equation}

We set $$W = \int_{\overline{\zeta}}^{\zeta}F(z, a, b)\,dz,$$ where $$F(z, a, b)=\frac
{i}{z^2-1}\,\log  \frac{4(z^2 + a^2)(z^2 + b^2)}{(1+a^2)(1+b^2)(z-z^2)^2}.$$ Now we
show that $W$ satisfies conditions (\ref{twelve}) and (\ref{thirteen}). So $W=V$ and
the theorem follows.

By the Leibniz formula  we have
\begin{equation} \label{fourteen}
\frac{\partial W}{\partial \alpha}=F(\zeta, a,
b)\frac{\partial \zeta}{\partial \alpha}-F(\overline{\zeta}, a,b)\frac{\partial
\overline{\zeta}}{\partial \alpha}+\, \int_{\overline{\zeta}}^{\zeta} \frac {\partial
F(z, a, b)}{\partial a}\, \frac{\partial a}{
\partial {\alpha}} \,dz
\end{equation}

We note that $F(\overline{\zeta}, a, b)=F(\zeta, a, b)=0$ if $\overline{\zeta} ,
\zeta, a$ and $b$ are the same as in the statement of the theorem. Moreover, since
$\alpha=2\,{\rm arccot\,}a$ we have $\displaystyle{\frac{\partial a}{\partial
{\alpha}}=-\frac{1+a^2}2}$ and

$$ \frac {\partial F(z, a, b)}{\partial a}\, \frac{\partial a}{
\partial {\alpha}} =-\frac{ia}{z^2+a^2}.
$$

Hence, by  Proposition \ref{ultima}, we obtain from (\ref{fourteen}) $$ \frac{\partial
W}{\partial \alpha}=-ia\,\int_{\overline{\zeta}}^{\zeta} \frac{\,dz}{z^2+a^2}=-
i\,\arctan\frac{a}{\zeta}+ i\,\arctan\frac{a}{\overline{\zeta}}=-\frac{r_{\alpha}}2.
$$

The equation $\displaystyle{\frac{\partial W}{\partial
\beta}=-\frac{r_{\beta}}2}$ can be obtained in the same way. The
boundary condition (\ref{thirteen}) for the function $W$ follows
from the integral formula.
\end{proof}
\medskip

\subsection{Particular cases and examples}

\begin{enumerate}
\item Case $\a=\b$.
In this case, Equation (\ref{prova}) splits into two quadratic
equations: $$(1+a^2)(z-z^2)+2(z^2+a^2)=0$$ and
$$(1+a^2)(z-z^2)-2(z^2+a^2)=0.$$ The first has two real roots
$z=-1$ and $z=2a^2/(a^2-1)$. The second has two non-real roots
$$z_{1,2}=\frac{1+a^2\pm\sqrt{1-22a^2-7a^4}}{2(3+a^2)}.$$ By
\cite{HLM4}, $\Delta=1-22a^2-7a^4$ is $<0$ in the hyperbolic case,
$=0$ in the Euclidean case and $>0$ in the spherical case. In the
Euclidean case we obtain $a^2={\rm
cot}^2(\a_0/2)=(\sqrt{128}-11)/7= 0.0448...$ and $a=a_0={\rm
cot\,}(\a_0/2)=0.2116...$ . So the cone-manifold is hyperbolic for
$0\le\a<\a_0=2.7243...$ and is Euclidean for $\a=\a_0$.

{}From (\ref{more}) we have $${\rm
Vol\,}W_2(\a,\a)=i\int_{z_1}^{z_2}\log\left[\frac{2(z^2+a^2)}{(z-z^2)(1+a^2)}\right]^2\frac{dz}{z^2-1}.$$
By differentiation with respect to $a$ and then by integration
with respect to $z$ we obtain $${\rm
Vol\,}W_2(\a,\a)=4\int_{a_0}^{a}{\rm
arctanh\,}\frac{\sqrt{7t^4+22t^2-1}}{t(5+t^2)}\frac{dt}{t^2+1}.$$

Since the integrand is pure imaginary for $0\le t<a_0$ we are able to compute the
volume in a more convenient way $${\rm Vol\,}W_2(\a,\a)=4\Re\int_{0}^{a}{\rm
arctanh\,}\frac{\sqrt{7t^4+22t^2-1}}{t(5+t^2)}\frac{dt}{t^2+1},$$
 where $a={\rm cot\,}\frac{\alpha}{2}.$

\item Case $\a=\b=\pi/2$. In this case equation (\ref{prova}) becomes $$(z+1)(z^2-z+2)=0.$$
Hence, the non-real roots are $$z_{1,2}=\frac{1\pm i\sqrt{7}}{2}$$
and $${\rm
Vol\,}W_2(\pi/2,\pi/2)=2i\int_{\frac{1-i\sqrt{7}}{4}}^{\frac{1+i\sqrt{7}}{4}}\log\frac{z^2+1}{z-z^2}\frac{dz}{z^2-1}=2.6667...$$

\item  Case $\a=\b=0$. Recall that $W_2(0,0)$ is the complete
hyperbolic manifold ${S}^3 \smallsetminus W_2$. By arguments
similar to the previous case, we obtain $${\rm
Vol\,}W_2(0,0)=2i\int_{\frac{1-i\sqrt{7}}{2}}^{\frac{1+i\sqrt{7}}{2}}\log\frac{2}{z-z^2}\frac{dz}{z^2-1}=5.3334...$$
Note that ${\rm Vol\,}W_2(0,0)=2\,{\rm Vol\,}W_2(\pi/2,\pi/2)$.

\item Case $\a=0$, $\b=\pi/3$. In this case equation (\ref{prova}) reduces to
$$(1+z)(3-3z+3z^2-z^3)=0.$$ Hence, the non-real roots are
$$z_{1,2}=1-\frac{1\pm i\sqrt{3}}{\sqrt[3]{4}}$$ and $${\rm
Vol\,}W_2(0,\pi/3)=i\int_{1-\frac{1+i\sqrt{3}}{\sqrt[3]{4}}}^{1-\frac{1-i\sqrt{3}}{\sqrt[3]{4}}}\log\frac{z^2+3}{(z-z^2)^2}\frac{dz}{z^2-1}=4.6165...$$
\end{enumerate}

The results of the above numerical calculation coincide with the
correspondent results obtained by Weeks's SnapPea program
\cite{We}.

\bigskip

\noindent{\bf Acknowledgement}

\noindent Work performed under the auspices of the G.N.S.A.G.A. of
I.N.d.A.M. (Italy) and the University of Bologna, funds for
selected research topics.

 \noindent The authors are thankful to
referees for helpful remarks and suggestions.

\bigskip\bigskip
%\newpage

\vspace{15 pt} {DMITRIY DEREVNIN, Novosibirsk State University, Novosibirsk, 630090,
Russia. E-mail: derevnin@mail.ru}

\vspace{15 pt} {ALEXANDER D. MEDNYKH, Sobolev Institute of
Mathematics, Novosibirsk, 630090, Russia. E-mail:
mednykh@math.nsc.ru}

\vspace{15 pt} {MICHELE MULAZZANI, Department of Mathematics and
C.I.R.A.M., University of Bologna, ITALY. E-mail:
mulazza@dm.unibo.it}

\end{document}